\documentclass[10pt]{article}
\font\smallit=cmti10

\usepackage{amssymb,latexsym,amsmath,epsfig,amsthm} 

\makeatletter

\renewcommand\section{\@startsection {section}{1}{\z@}
{-30pt \@plus -1ex \@minus -.2ex}
{2.3ex \@plus.2ex}
{\normalfont\normalsize\bfseries}}

\renewcommand\subsection{\@startsection{subsection}{2}{\z@}
{-3.25ex\@plus -1ex \@minus -.2ex}
{1.5ex \@plus .2ex}
{\normalfont\normalsize\bfseries}}

\renewcommand{\@seccntformat}[1]{\csname the#1\endcsname. }

\makeatother

\newtheorem*{theorem*}{Theorem}
\newtheorem{lemma}{Lemma}

\newtheorem*{corollary*}{Corollary}
\theoremstyle{remark}
\newtheorem*{remark*}{Remark}
\begin{document}
\begin{center}
\uppercase{\bf An infinite collection of quartic polynomials whose products of consecutive values are not perfect squares}
\vskip 20pt
{\bf Konstantinos Gaitanas}\\
{\smallit Department of  Applied Mathematical and Physical Sciences, National Technical University of Athens, Greece}\\
{\tt raffako@hotmail.com}\\
\vskip 10pt
\end{center}
\vskip 30pt
\vskip 30pt

\centerline{\bf Abstract}
\noindent
Using an elementary identity, we prove that for infinitely many polynomials $P(x)\in \mathbb{Z}[X]$ of fourth degree, the equation $\prod\limits_{k=1}^{n}P(k)=y^2$ has finitely many solutions in $\mathbb{Z}$. We also give an example of a quartic polynomial for which the product of it's first consecutive values is infinitely often a perfect square.
\pagestyle{myheadings}  
\thispagestyle{empty} 
\baselineskip=12.875pt 
\vskip 30pt
\section{Introduction} 
Over the last few years, there has been a growing interest in identifying if certain product sequences contain perfect squares. In 2008 Javier Cilleruelo \cite{1} proved that the product $(1^2+1)(2^2+1)\cdots (n^2+1)$ is a square only for $n=3$. Soon after, Jin-Hui Fang \cite{2} achieved to prove that both of the products $\prod\limits_{k=1}^{n}(4k^2+1)$ and $\prod\limits_{k=1}^{n}\big(2k(k-1)+1\big)$ are never squares. There are not many similar results for quadratic polynomials. However, in a recent paper \cite{3} two certain cases of quartic polynomials were settled. In this paper we will prove using elementary arguments that there is actually an infinite collection of quartic polynomials $P(x)$ such that the product $\displaystyle\prod\limits_{k=1}^{n}P(k)$ is a square finitely often. At the end of the paper we discuss some cases that can be handled by this method. We begin with a polynomial identity which is the key ingredient throughout this article.

\begin{lemma}
Let $f(x)=x^2+ax+ b$ be a quadratic polynomial. For every $x\in \mathbb{R}$ the following formula is valid:
$$f\big(f(x) + x\big) = f(x)f(x + 1)$$.
\begin{proof}
We can verify this just by doing elementary manipulations but we will prove the lemma using a clever observation. Since $f(x)$ is a polynomial of second degree, Taylor's formula gives $f\big(f(x) + x\big)=f(x)+\frac{f'(x)f(x)}{1!}+f^2(x)$. This is equal to $f(x)\big(1+f'(x)+f(x) \big)$. But $1+f'(x)+f(x)=1+2x+a+x^2+ax+b=(x+1)^2+a(x+1)+b=f(x+1)$. Hence we have: $f\big(f(x) + x\big) = f(x)f(x + 1)$.
\end{proof}
\end{lemma}
This simple formula will play a key role in the proof of the main theorem. For convenience of notation we set $f\big(f(k) + k\big)=P(k)=f(k)f(k+1)$. It can be seen that $P(k)=k^4+2(a+1)k^3+\big((a+1)^2+2b+a)\big)k^2+(a+1)(2b+a)k+b^2+ab+b$. In the proof of the main theorem, we require $a$ and $b$ to obey a certain restriction. Under this restriction, we are able to prove that equation (1) has finitely many solutions.
\section{Main Results}
\begin{theorem*}
Let $a, b, m\in \mathbb{Z}$ and $a+b+1=m^2$. Then the diophantine equation \begin{equation}\displaystyle\prod\limits_{k=1}^{n} P(k)=y^2\label{1}\end{equation} has finitely many solutions.
\begin{proof}
Using lemma 1 we can rewrite equation ~(\ref{1}) as $f(1)f(2)f(2)f(3)\cdots f(n)f(n+1)=y^2$ which reduces to $f(1)f(n+1)\displaystyle\prod\limits_{k=2}^{n}(f(k))^2=y^2$. Since $f(1)=a+b+1=m^2$ we conclude that $f(n+1)=\frac{y^2 }{m^2\prod\limits_{k=2}^{n}(f(k))^2}$. It becomes clear that equation ~(\ref{1}) is satisfied whenever $f(n+1)$ is a perfect square. It remains to prove that among the values of $f(k)$ occur finitely many squares. Write\begin{equation} k^2+ak+b=z^2 \label{2}\end{equation} for some $z\in \mathbb{Z}$. This means that for sufficiently large $k$, $k^2<z^2<(k+2a)^2$ if $a> 0$ or, $(k+2a)^2<z^2<k^2$ if $a< 0$. (If $a=0$ then equation ~(\ref{2}) transposes to $(z-x)(z+x)=b$ which clearly has finitely many solutions). Both of the inequalities yield $z=k+c$ for some $c\in \mathbb{Z}$ with $|c|<|2a|$. So, ~(\ref{2}) becomes $k^2+ak+b=(k+c)^2$ which has finitely many solutions as the reader may easily verify.
\end{proof}
\end{theorem*}
It suffices to choose some nice values for a and b in order to demonstrate the theorem. Choosing $(a, b) = (-1, 1)$ we have $f(k) = k^2-k + 1$ hence the following: 
\begin{corollary*}
$\displaystyle\prod\limits_{k=1}^{n}(k^4+k^2+1)$ is a square only for $n=1$.
\begin{proof}
If $(a, b)=(-1, 1)$ then $f(1)=1^2$. Repeating the previous arguments, it suffices to show that $k^2-k+1=y^2$ has one solution. Indeed, if $k^2-k+1=y^2$ then we must have $k^2\le y^2<(k+1)^2$ which yields $y=k$ and so $k=1$. The claim follows.
\end{proof}
\end{corollary*}
\begin{remark*}
Arguing as in the previous section, we may present an example which shows that equation ~(\ref{1}) has infinitely many solutions. Choosing $(a, b)=(-4, 2)$ we have $f(k)=k^2-4k+2$ and $P(k)=\big(k(k-3)\big)^2-2$. We can prove that the product $\displaystyle\prod\limits_{k=4}^{n}\Big(\big(k(k-3)\big)^2-2\Big)$ is a square infinitely often. Here we start with $k=3$ to omit any trivial case in which the product has negative factors. The product is a square if $f(4)f(n+1)=2(n-1)^2-4=y^2$. It is a routine matter to prove that both $y$ and $n-1$ must be even. Thus, equation can be written as $(\frac{y}{2})^2-2(\frac{n-1}{2})^2=-1$ which is a special case of the negative Pell equation $X^2-2Y^2=-1$. This equation has the fundamental solution $(1, 1)$ and all it's positive solutions can be found by taking odd powers of $1+\sqrt 2$. The positive solutions are $(X_n,Y_n)$ where $X_n+Y_n\sqrt 2=(1+\sqrt 2)^{2n-1}$. The next solution is $(X_2, Y_2)=(7, 5)$ which gives $n=11$. As an example we can verify that $\displaystyle\prod\limits_{k=4}^{11}\Big(\big(k(k-3)\big)^2-2\Big)=246988938224^2$
\end{remark*}
\makeatletter
\renewcommand{\@biblabel}[1]{[#1]\hfill}
\makeatother

\end{document}